\documentclass[12pt,twoside]{amsart}
\usepackage{amssymb}
\usepackage{graphicx}
\usepackage{amsmath}

\newtheorem{theorem}{Theorem}[section]
\newtheorem{corollary}[theorem]{Corollary}
\newtheorem{definition}[theorem]{Definition}

\newtheorem{lemma}[theorem]{Lemma}
\newtheorem{notation}[theorem]{Notation}
\newtheorem{proposition}[theorem]{Proposition}

\begin{document}

\begin{center}
{\bf\Large On the order of a non-abelian representation
group of a slim dense near hexagon}\\

\vskip 3cm

{\bf\large Binod Kumar Sahoo}\footnote{Supported by DAE Grant
39/3/2000-R\&D-II (NBHM Fellowship), Govt. of India.}\\
(binodkumar@gmail.com)\smallskip\\ {\it and}\smallskip\\
{\bf\large N. S. Narasimha Sastry}\\
(nnsastry@gmail.com)\bigskip\\
Statistics \& Mathematics Unit\\
Indian Statistical Institute\\
8th Mile, Mysore Road\\
R.V. College Post\\
Bangalore-560059, India\smallskip\\
(21 April 2006)
\end{center}
 \vskip 2cm

\begin{quote}
{\bf Abstract.} We show that, if the representation group $R$ of a
slim dense near hexagon $S$ is non-abelian, then $R$ is of exponent
4 and $|R|=2^{\beta}$, $1+NPdim(S)\leq \beta\leq 1+dimV(S)$, where
$NPdim(S)$ is the near polygon embedding dimension of $S$ and
$dimV(S)$ is the dimension of the universal representation module
$V(S)$ of $S$. Further, if $\beta =1+NPdim(S)$, then $R$ is an
extraspecial 2-group (Theorem \ref{main}).

\vskip 2cm

\smallskip

\end{quote}

\pagebreak

\title[Order of a non-abelian representation group]{}
\author[B. K. Sahoo]{}
\address[]{}
\author[N. S. N. Sastry]{}
\address[]{}
\maketitle

\section{Introduction}

A \textit{partial linear space} is a pair $S=(P,L)$ consisting of a
nonempty `point-set' $P$ and a nonempty `line-set' $L$ of subsets of
$P$ of size at least 3 such that any two distinct points $x$ and $y$
are in at most one line. Such a line, if it exists, is written as
$xy$, $x$ and $y$ are said to be \textit{collinear} and written as
$x\sim y$. If $x$ and $y$ are not collinear, we write $x\nsim y$. If
each line contains exactly three points, then $S$ is \textit{slim}.
For $ x\in P$ and $A\subseteq P$, we define $x^{\perp
}=\left\{x\right\}\cup \left\{y\in P:x\sim y\right\}$ and $A^{\perp
}=\underset{x\in A}{\cap }x^{\perp }$. If $P^{\perp }$ is empty,
then $S$ is \textit{non-degenerate}. A subset of $P$ is a
\textit{subspace} of $S$ if any line containing at least two of its
points is contained in it. For a subset $X$ of $P$, the
\textit{subspace} $ \left\langle X\right\rangle$ \textit{generated}
by $X$ is the intersection of all subspaces of $S$ containing $X$. A
\textit{geometric hyperplane} of $S$ is a subspace of $S$, different
from the empty set and $P$, that meets every line nontrivially. The
graph $\Gamma \left( P\right)$ with vertex set $P$, two distinct
points being {\it adjacent} if they are collinear in $S$, is the
\textit{collinearity graph} of $S$. For $x\in P$ and an integer $i$,
we write
\begin{eqnarray*}
\Gamma _{i}(x) &=&\{ y\in P:d(x,y) =i\},\\
\Gamma _{\leq i}(x) &=&\{y\in P:d(x,y) \leq i\},
\end{eqnarray*}
where $d(x,y)$ denotes the \textit{distance} between $x$ and $y$ in
$\Gamma (P)$. The \textit{diameter} of $S$ is the diameter of
$\Gamma \left( P\right)$. If $\Gamma \left( P\right) $ is connected,
then $S$ is a \textit{connected} point-line geometry.

\subsection{Representation of a partial linear space}

Let $S=(P,L)$ be a connected slim partial linear space. If $x,y\in
P$ and $x\sim y$, we define $x*y$ by $xy =\{x,y,x*y\}$.

\begin{definition}
\label{definition}(\cite{IPS}, p.525) A representation $(R,\psi)$ of
$S$ with representation group $R$ is a mapping $\psi:x\mapsto
\langle r_{x}\rangle$ from $P$ into the set of subgroups of order
$2$ of $R$ such that the following hold:
\begin{enumerate}
\item[$(i)$] $R$ is generated by $Im(\psi)$.

\item[$(ii)$] If
$l=\{ x,y,x*y\} \in L$, then $\{1,r_{x},r_{y},r_{x*y}\}$ is a Klein
four group.
\end{enumerate}
\end{definition}

We set $R_{\psi}=\{ r_{x}:x\in P\}$. The representation $(R,\psi)$
is \textit{faithful} if $\psi $ is injective. A representation
$(R,\psi)$ of $S$ is \textit{abelian} or \textit{non-abelian}
according as $R$ is abelian or not. Note that, in \cite{IPS},
`non-abelian representation' means `the representation group is not
necessarily abelian'.

For an abelian representation, the representation group can be
considered as vector space over $F_{2}$, the field with two
elements. For each connected slim partial linear space $S$, there
exists a unique abelian representation $\rho _{0}$ of $S$ such that
any other abelian representation of $S$ is a composition of $\rho
_{0}$ and a linear mapping (see \cite{R}). $\rho _{0}$ is called the
{\it universal abelian representation} of $S$. The $F_{2}$ vector
space $V(S)$ underlying the universal abelian representation is
called the \textit{universal representation module} of $S$.
Considering $V(S)$ as an abstract group with the group operation
$+$, it has the presentation
\begin{equation*}
\begin{array}{ccc}
V\left( S\right)  & = & \langle v_{x}:x\in P;\text{ }2v_{x}=0;
\text{ }v_{x}+v_{y}=v_{y}+v_{x}
\text{ for }x,y\in P; \\
&  & \text{and }v_{x}+v_{y}+v_{x\ast y}=0\text{ if }x\sim y\rangle
\end{array}
\end{equation*}
and $\rho_{0}$ is defined by $\rho_{0}(x)= v_{x}$ for $x\in P$.

A representation $(R_{1},\psi _{1})$ of $S$ is a \textit{cover} of
a representation $(R_{2},\psi _{2})$ of $S$ if there exist an
automorphism $\beta$ of $S$ and a group homomorphism $\varphi
:R_{1}\rightarrow R_{2}$ such that $\psi _{2}(\beta( x)) =\varphi
(\psi _{1}(x))$ for every $x\in P$. If $S$ admits a non-abelian
representation, then there is a \textit{universal representation}
$(R(S),\psi_{S})$ which is the cover of every other representation
of $S$. The universal representation is unique (see \cite{I}, p.
306) and the \textit{universal representation group} $R\left(
S\right) $ of $S $ has the presentation:
\begin{equation*}
R\left( S\right) =\langle r_{x}:x\in
P,r_{x}^{2}=1,r_{x}r_{y}r_{z}=1\text{ if }\left\{ x,y,z\right\}
\in L\rangle .
\end{equation*}

Whenever we have a representation of $S$, the group spanned by the
images of the points is a quotient of $R(S)$. Further,

\begin{lemma}
\label{order-calculation}(\cite{IPS}, p.525) $V(S) =R(S) /[R(S)
,R(S) ]$.
\end{lemma}

The general notion of a representation group of a finite partial
linear space with $p+1$ points per line for a prime $p$ was
introduced by Ivanov \cite{I} in his investigations of Petersen
and Tilde geometries (motivated in large measure by questions
about the Monster and Baby Monster finite simple groups). A
sufficient condition on the partial linear space and on the
non-abelian representation of it is given in \cite{SS} to ensure
that the representation group is a finite $p$-group. For a
detailed survey on non-abelian representations, we refer to
\cite{I}, also see (\cite{SS}, Sections 1 and 2). We study the
possible orders of a non-abelian representation group of a slim
dense near near hexagon here (Theorem \ref{main}).

\subsection{Near $2n$-gons}

A \textit{near $2n$-gon} is a connected non-degenerate partial
linear space $S=(P,L)$ of diameter $n$ such that for each point-line
pair $(x,l)\in P\times L$, $x$ is nearest to exactly one point of
$l$. Near 4-gons are precisely \textit{generalized quadrangles} (GQ,
for short); that is, non-degenerate partial linear spaces such that
for each point-line pair $(x,l)$, $x\notin l$, $x$ is collinear with
exactly one point of $l$.

Let $S=(P,L)$ be a near $2n$-gon. Then the sets $S(x)=\Gamma_{\leq
n-1}(x)$, $x\in P$, are \textit{special} geometric hyperplanes. A
subset $C$ of $P$ is \textit{convex} if every shortest path in
$\Gamma (P)$ between two points of $C$ is entirely contained in $C$.
A \textit{quad} is a convex subset of $P$ of diameter $2$ such that
no point of it is adjacent to all other points of it. If
$x_{1},x_{2}\in P$ with $d(x_{1},x_{2})=2$ and
$|\{x_{1},x_{2}\}^{\perp}|\geq 2$, then $x_{1}$ and $x_{2}$ are
contained in a unique quad, denoted by $Q(x_{1},x_{2})$, which is a
generalized quadrangle (\cite{SY}, Proposition 2.5, p.10). Thus, a
quad is a subspace.

A near $2n$-gon is called \textit{dense} if every pair of points at
distance 2 are contained in a quad. In a dense near $2n$-gon, the
number of lines through a point is independent of the point
(\cite{BW}, Lemma 19, p.152). We denote this number by $t+1$. A near
$2n$-gon is said to have \textit{parameters $(s,t)$} if each line
contains $s+1$ points and each point is contained in $t+1$ lines. A
near 4-gon with parameters $(s,t)$ is written as $(s,t)$-GQ.

\begin{theorem}
\label{SY1}(\cite{SY}, Proposition 2.6, p.12) Let $S=(P,L)$ be a
near $2n$-gon and $Q$ be a quad in $S$. Then, for $x\in P$, either
\begin{enumerate}
\item[$(i)$] there is a unique point $y\in Q$ closest to $x$
(depending on $x$) and $d(x,z)=d(x,y)+d(y,z)$ for all $z\in Q$; or
\item[$(ii)$] the points in $Q$ closest to $x$ form an ovoid
$\mathcal{O}_{x}$ of $Q$.
\end{enumerate}
\end{theorem}

The point-quad pair $(x,Q)$ in Theorem \ref{SY1} is called
\textit{classical} in the first case and \textit{ovoidal} in the
later case. A quad $Q$ is \textit{classical} if $(x,Q)$ is classical
for each $x\in P$, otherwise it is \textit{ovoidal}.

\subsection{Slim dense near hexagons}

A near 6-gon is called a \textit{near hexagon}. Let $S=(P,L)$ be a
slim dense near hexagon. For $x,y\in P$ with $d(x,y)=2$, we write
$|\Gamma_{1}(x)\cap \Gamma_{1}(y)|$ as $t_{2}+1$ (though this
depends on $x,y$). We have, $t_{2}<t$. A quad in $S$ is \textit{big}
if it is classical. Thus, if $Q$ is a big quad in $S$, then each
point of $S$ has distance at most one to $Q$. We say that a quad $Q$
is of \textit{type} $(2,t_{2})$ if it is a $(2,t_{2})$-GQ.

\begin{theorem}\label{classification-result}
(\cite{BCHW}, Theorem 1.1, p.349) Let $S=(P,L)$ be a slim dense near
hexagon. Then $S$ is necessarily finite and is isomorphic to one of
the eleven near hexagons with parameters as given below.
\begin{center}
$
\begin{array}{|r|c|c|c|c|c|c|c|c|} \hline
& |P|  & t & t_{2} &  dim V(S) & NPdim(S) & a_{1} & a_{2} & a_{4} \\
\hline
(i) & 759 & 14 & 2 & 23 & 22 & - & 35 & - \\ \hline
(ii) &
729 & 11 & 1 & 24 & 24 & 66 & - & - \\ \hline
(iii)  & 891 & 20 &
4^{\star} & 22 & 20 & - & - & 21 \\ \hline
(iv)  & 567 & 14 &
2,4^{\star} & 21 & 20 & - & 15 & 6 \\ \hline
(v)  & 405 & 11 &
1,2,4^{\star} & 20 & 20 & 9 & 9 & 3 \\ \hline
(vi)  & 243 & 8 &
1,4^{\star} & 18 & 18 & 16 & - & 2 \\ \hline
(vii)  & 81 & 5 &
1,4^{\star} & 12 & 12 & 5 & - & 1
\\ \hline
(viii)  & 135 & 6 & 2^{\star} & 15 & 8 & - & 7 & - \\
\hline (ix)  & 105 & 5 & 1,2^{\star} & 14 & 8 & 3 & 4 & - \\ \hline
(x)  & 45 & 3 & 1,2^{\star} & 10 & 8 & 3 & 1 & - \\ \hline
(xi)  &
27 & 2 & 1^{\star} & 8 & 8 & 3 & - & - \\ \hline
\end{array}
$
\end{center}
\end{theorem}

Here, $NPdim(S)$ is the $F_{2}$-rank of the matrix $A_{3}:P\times
P\longrightarrow \{0,1\}$ defined by $A_{3}(x,y)=1$ if $d(x,y)=3$
and zero otherwise. We add a star if and only if the corresponding
quads are big. The number of quads of type $(2,r)$, $r=1,2,4$,
containing a point of $S$ in indicated by $a_{r}$. A `--' in a
column means that $a_{r}=0$.

For a description of the near hexagons $(i)-(iii)$ see \cite{SY} and
for $(iv)-(xi)$ see \cite{BCHW}. However, the parameters of these
near hexagons suffice for our purposes here. For other
classification results about slim dense near polygons, see
\cite{Bart-book} and \cite{V}.

\subsection{Extraspecial 2-groups}

A finite $2$-group $G$ is \textit{extraspecial} if its Frattini
subgroup $ \Phi \left( G\right) ,$ the commutator subgroup
$G^{\prime }$ and the center $Z\left( G\right) $ coincide and have
order $2$.

An extraspecial $2$-group is of exponent 4 and order $2^{1+2m}$ for
some integer $m\geq 1$ and the maximum of the orders of its abelian
subgroups is $2^{m+1}$ (see \cite{DH}, section 20, p.78,79). An
extraspecial 2-group $G$ of order $2^{1+2m}$ is a central product of
either $m$ copies of the dihedral group $D_{8}$ of order 8 or $m-1$
copies of $D_{8}$ with a copy of the quaternion group $Q_{8}$ of
order 8. In the first case, $G$ possesses a maximal elementary
abelian subgroup of order $2^{1+m}$ and we write $G=2_{+}^{1+2m}$.
If the later holds, then all maximal abelian subgroups of $G$ are of
the type $2^{m-1}\times 4$ and we write $G=2_{-}^{1+2m}$.

\begin{notation}
For a group $G$, $G^{*}=G\setminus \{1\}$.
\end{notation}

\subsection{The main result}

In this paper, we prove

\begin{theorem}
\label{main}Let $S=(P,L)$ be a slim dense near hexagon and
$(R,\psi)$ be a non-abelian representation of $S$. Then
\begin{enumerate}
\item[$(i)$] $R$ is a finite 2-group of exponent 4 and order
$2^{\beta}$, where $1+NPdim (S)\leq \beta\leq 1+dim V(S)$.

\item[$(ii)$] If $\beta = 1+NPdim (S)$, then $R$ is an
extraspecial 2-group. Further, $R=2_{+}^{1+NPdim(S)}$ except for the
near hexagon $(vi)$ in Theorem \ref{classification-result}. In that
case, $R=2_{-}^{1+NPdim(S)}$.
\end{enumerate}
\end{theorem}

Section 2 is about slim dense near hexagons. In Section 3, we
study representations of $(2,t)$-GQs. In Section 4, we study the
non-abelian representations of slim dense near hexagons. In
section 5 we prove Theorem \ref{main}.

\section{Elementary Properties}

Let $S=(P,L)$ be a slim dense near hexagon. Since a (2,4)-GQ admits
no ovoids, every quad in $S$ of type $(2,4)$ is big (see Theorem
\ref{SY1}).

\begin{lemma}\label{elementary1}
(\cite{BCHW}, p.359) Let $Q$ be a quad in $S$ of type $(2,t_{2})$.
Then $|P|\geq |Q|(1+2(t-t_{2}))$. Equality holds if and only if $Q$
is big. In particular, if a quad in $S$ of type $(2,t_{2})$ is big
then so are all quads in $S$ of that type.
\end{lemma}

Let $Q_{1}$ and $Q_{2}$ be two disjoint big quads in $S$.

\begin{lemma}\label{isomorphism}
(\cite{BCHW}, Proposition 4.3, p.354) Let $\pi$ be the map from
$Q_{1}$ to $Q_{2}$ which takes $x$ to $z_{x}$, where $x\in Q_{1}$
and $z_{x}$ is the unique point in $Q_{2}$ at a distance one from
$x$. Then
\begin{enumerate}
\item[$(i)$] $\pi$ is an isomorphism from $Q_{1}$ to $Q_{2}$.

\item[$(ii)$] The set
$Q_{1}*Q_{2}=\{x*z_{x}:x\in Q_{1}\}$ is a big quad in $S$.
\end{enumerate}
\end{lemma}

Let $Y$ be the subspace of $S$ generated by $Q_{1}$ and $Q_{2}$.
Note that $Y$ is isomorphic to the near hexagon $(ix),(x)$ or
$(vii)$ according as $Q_{1}$ and $Q_{2}$ are GQs of type (2,1),
(2,2) or (2,4). Let $\{i,j\} =\{1,2\}$. For $x\in P\setminus Y$, we
denote by $x^{j}$ the unique point in $Q_{j}$ at a distance 1 from
$x$. For $y\in Q_{i}$, $z_{y}\in Q_{j}$ is defined as in Lemma
\ref{isomorphism}. The following elementary results are useful for
us.

\begin{proposition}
\label{done1(vi)}For $x\in P\setminus Y$, $d(z_{x^{i}},x^{j})=1$ and
$d(z_{x^{1}},z_{x^{2}})=d(x^{1},x^{2})=2$; that is,
$\{x^{1},z_{x^{1}},x^{2},z_{x^{2}}\}$ is a quadrangle in
$\Gamma(P)$.
\end{proposition}

\begin{proof}[\bf Proof]
Since $x\in \Gamma _{1}(x^{1}) \cap \Gamma _{1}(x^{2})$,
$d(x^{1},x^{2}) =2$. Further, $d( x^{i},x^{j}) =d(x^{i},z_{x^{i}})
+d( z_{x^{i}},x^{j})$. So $d( z_{x^{i}},x^{j}) =1$ and
$d(z_{x^{1}},z_{x^{2}}) =2$.
\end{proof}

\begin{proposition}\label{done3(vi)}
Let $l$ be a line of $S$ disjoint from $Y$ and $x,y\in l$, $x\neq
y$. Then, $x^{1}y^{1}=x^{1}z_{x^{2}}$ if and only if $
x^{2}y^{2}=x^{2}z_{x^{1}}$. In fact, if $x^{1}y^{1}=x^{1}z_{x^{2}}$,
then $(y^{1},y^{2}) =(z_{x^{2}},x^{2}\ast z_{x^{1}})$ or $(x^{1}\ast
z_{x^{2}},z_{x^{1}})$.
\end{proposition}

\begin{proof}[\bf Proof]
$ x^{j}y^{j}=x^{j}z_{x^{i}}$ if and only if
$y^{j}\in\{z_{x^{i}},x^{j}\ast z_{x^{i}}\}$. If $y^{j}=x^{j}\ast
z_{x^{i}}$, then $y^{i}\sim x^{i}\ast z_{x^{j}}$, because
$2=d(y^{j},y^{i})=d(y^{j},x^{i}\ast z_{x^{j}})+d(x^{i}\ast
z_{x^{j}},y^{i})$. Since $y^{i}\sim x^{i}$, it follows that $y^{i}$
is a point in the line $x^{i}z_{x^{j}}$ and $y^{i}=z_{x^{j}}$.

If $y^{j}=z_{x^{i}}$, then applying the above argument to $(x\ast
y)^{j}=x^{j}\ast z_{x^{i}}$, we get $(x\ast y)^{i}=z_{x^{j}}$ and
$y^{i}=x^{i}\ast z_{x^{j}}$.
\end{proof}

An immediate consequence of Proposition \ref{done3(vi)} is the
following.

\begin{corollary}\label{done(vi)-corollary}
Let $l$ be a line of $S$ disjoint from $Y$ and $x,y\in l$ with
$x\neq y$. Then, $d(z_{x^{i}},z_{y^{j}})\leq 2$ or 3 according as
the lines $x^{j}y^{j}$ and $x^{j}z_{x^{i}}$ coincide or not.
\end{corollary}

\begin{proposition}\label{done(vi)-proposition}
Let $Q$ be a big quad in $S$ disjoint from $Y$. For $x,y\in Q$ with
$x\nsim y$, $(d(z_{x^{1}},z_{y^{2}}),d(z_{x^{2}},z_{y^{1}}))=(2,3)$
or (3,2).
\end{proposition}

\begin{proof}[\bf Proof]
By Lemma \ref{isomorphism}, there exist $w\in \{x,y\} ^{\perp }$ in
$Q$ such that $x^{1}w^{1}=x^{1}z_{x^{2}}$. By Proposition
\ref{done3(vi)}, $(w^{1},w^{2}) =(z_{x^{2}},x^{2}\ast z_{x^{1}})$ or
$(x^{1}\ast z_{x^{2}},z_{x^{1}})$. Assume that $(w^{1},w^{2})
=(z_{x^{2}},x^{2}\ast z_{x^{1}})$. Then,
$d(z_{x^{2}},z_{y^{1}})=d(w^{1},z_{y^{1}})=d(w^{1},z_{w^{1}})+d(z_{w^{1}},z_{y^{1}})=2$.
Now, $y^{2}\sim w^{2}$ and $y^{2}\nsim x^{2}$ in $Q_{2}$ implies
that $x^{1}\nsim z_{y^{2}}$. So $d(x^{1},z_{y^{2}})=2$ and
$d(z_{x^{1}},z_{y^{2}})=d(z_{x^{1}},x^{1})+d(x^{1},z_{y^{2}})=3$. A
similar argument holds if $(w^{1},w^{2})=(x^{1}\ast
z_{x^{2}},z_{x^{1}})$.
\end{proof}

\section{Representations of $(2,t)$-GQs}

Let $S=\left( P,L\right) $ be a $\left( 2,t\right) $-GQ. Then $P$ is
finite and $t=1,2$ or $4$. For each value of $t$ there exists a
unique generalized quadrangle, up to isomorphism (\cite{C}, Theorem
7.3, p.99). A $k$\textit{-arc} of $S$ is a set of $k$ pair-wise
non-collinear points of $S$. A $k$-arc is \textit{complete} if it is
not contained in a $\left( k+1\right) $-arc. A point $x$ is a
\textit{center} of a $k$-arc if $x$ is collinear with every point of
it. An \textit{ovoid} of $S$ is a $k$-arc meeting each line of $S$
non-trivially. A \textit{spread} of $S$ is a set $K$ of lines of $S$
such that each point of $S$ is in a unique member of $K$. If $O$
(resp., $K$) is an ovoid (resp., spread) of $S$, then $|O|=1+2t$
(resp., $|K|=1+2t$).

Since each line contains three points, each pair of non-collinear
points of $S$ is contained in a $(2,1)$-subGQ of $S$. For $t=1,2$, a
$(2,t)$-subGQ of $S$ and a point outside it generate a
$(2,2t)$-subGQ in $S$. Minimum number of generators of a $(2,t)$-GQ
is 4 if $t=1$, 5 if $t=2$ and 6 if $t=4$.

\subsection{$(2,2)$-GQ}

Let $S=(P,L)$ be a $(2,2)$-GQ. For any 3-arc $T$ of $S$,
$|T^{\perp}| =1$ or $3$. Further, $|T^{\perp }| =1$ if and only if
$T$ is contained in a unique $(2,1)$-subGQ of $S$; and $|T^{\perp}|
=3$ if and only if $T$ is a complete $3$-arc. If $S$ admits a
$k$-arc, then $k\leq 5$. $S$ contains six 5-arcs (that is, ovoids).
Each ovoid is determined by any two of its points. Each point of $S$
is in two ovoids and the intersection of two distinct ovoids is a
singleton. Any two non-collinear points of $S$ are in a unique ovoid
of $S$ and also in a unique complete 3-arc of $S$. Any incomplete
3-arc of $S$ is contained in a unique ovoid. Any $4$-arc of $S$ is
not complete and is contained in a unique ovoid. The intersection of
two distinct complete 3-arcs of $S$ is empty or a singleton.

\label{model1}{\bf A model for the $(2,2)$-GQ:} Let $\Omega
=\{1,2,3,4,5,6\}$. A \textit{factor} of $\Omega $ is a set of three
pair-wise disjoint 2-subsets of $\Omega$. Let $\mathcal{E}$ be the
set of all 2-subsets of $\Omega $ and $\mathcal{F}$ be the set of
all factors of $\Omega$. Then $|\mathcal{E}|=|\mathcal{F}|=15$ and
the pair $(\mathcal{E},\mathcal{F})$ is a $(2,2)$-GQ.

\subsection{$(2,4)$-GQ}

Let $S=(P,L)$ be a $(2,4)$-GQ. Each 3-arc of $S$ has three centers
and is contained in a unique $(2,1)$-subGQ of $S$. So any $4$-arc of
$S$ is contained in a unique $(2,2)$-subGQ of $S$. If $S$ admits a
$k$-arc, then $0\leq k\leq 6$. So $S$ has no ovoids. $S$ admits two
disjoint $6$-arcs. A $5$-arc of $S$ is complete if and only if it is
contained in a unique $(2,2)$-subGQ of $S$. Each incomplete 5-arc
has exactly one center and each complete $5$-arc of $S$ has exactly
two centers. Each $4$-arc has two centers and is contained in a
unique complete $5$-arc and in a unique complete $6$-arc. Each 3-arc
of $S$ has 3 centers and is contained in a unique (2,1)-subGQ of
$S$.

{\bf A model for the $(2,4)$-GQ:} Let $\Omega$, $\mathcal{E}$ and
$\mathcal{F}$ be as in the model of a (2,2)-GQ. Let $\Omega ^{\prime
}=\{1',2',3',4',5',6'\}$. Take
$$P =\mathcal{E}\cup \Omega \cup \Omega ^{\prime };\text{ } L =\mathcal{F}\cup
\{\{i,\{i,j\},j'\}:1\leq i\neq j\leq 6\}.$$ Then $|P|=27$, $|L|=45$
and the pair $(P,L)$ is a (2,4)-GQ.

\subsection{Representations}

Let $S=(P,L)$ be a $(2,t)$-GQ and $(R,\psi)$ be a representation of
$S$.

\begin{proposition}
\label{begin}$R$ is an elementary abelian 2-group.
\end{proposition}

\begin{proof}[\bf Proof]
Let $x,y\in P$ and $x\nsim y$. Let $T$ be a $(2,1)$-subGQ of $S$
containing $x$ and $y$. Let $\{x,y\}^{\perp }$ in $T$ be $\{a,b\}$.
Then $[r_{x},r_{y}]=1$, because $r_{b}r_{y}=r_{y}r_{b}$,
$r_{b}r_{x}=r_{x}r_{b}$ and $r_{(a\ast x) \ast (b\ast y)}=r_{(a\ast
y) \ast ( b\ast x)}$. So $R$ is abelian.
\end{proof}

For the rest of this section we assume that $\psi $ is faithful.

\begin{proposition}
\label{order for GQ}The following hold:
\begin{enumerate}
\item[$(i)$] $\left| R\right| =2^{4}$ if $t=1$;
\item[$(ii)$]
$\left| R\right| =2^{4}$ or $2^{5}$ if $t=2$, and both
possibilities occur;
\item[$(iii)$] $ \left| R\right| =2^{6}$ if
$t=4$.
\end{enumerate}
\end{proposition}

\begin{proof}[\bf Proof]
Since $S$ is generated by $k$ points, where $(t,k)\in
\{(1,4),(2,5),(4,6)\}$, $F_2$-dimension of $R$ is at most $k$. So
$|R|\leq 2^{k}$.

$(i)$ If $t=1$, then $| R| \geq 2^{4}$ because $|P| =9$ and $\psi$
is faithful. So $|R|=2^4$.

$(ii)$ If $t=2$, then $|R|\geq 2^4$ because $S$ contains a
$(2,1)$-subGQ. The rest follows from the fact that $S$ has a
symplectic embedding in a $F_{2}$-vector space of dimension 4 and as
well as an orthogonal embedding in a $F_{2}$-vector space of
dimension 5.

To prove $(iii)$ we need Proposition \ref{converse} below which is a
partial converse to the fact that if $x\sim y$, $x,y\in P$, then
$r_xr_y\in R_\psi $.
\end{proof}

\begin{proposition}
\label{converse}Assume that $(t,|R|) \neq (2,2^{4})$. If $r_xr_y\in
R_\psi$ for distinct $x,y\in P$, then $x\sim y$.
\end{proposition}

\begin{proof}[\bf Proof]
Let $z\in P$ be such that $r_z=r_xr_y$. If $x\nsim y$, then $
T=\{x,y,z\}$ is a 3-arc of $S$ because $\psi$ is faithful. There is
no (2,1)-subGQ of $S$ containing $T$ because the subgroup of $R$
generated by the image of such a GQ is of order $2^{4}$ (Proposition
\ref{order for GQ}$(i)$). Every 3-arc of a $ (2,4)$-GQ is contained
in a unique (2,1)-subGQ. So $t=2$ and $T$ is a complete 3-arc. Let
$Q$ be a (2,1)-subGQ of $S$ containing $x$ and $y$. Then $z\notin Q$
and $P=\langle Q,z\rangle$. Since $r_z\in \left\langle \psi (Q)
\right\rangle$, $
|R| =|\left\langle \psi (Q) \right\rangle| =2^{4},$ a contradiction to $%
\left( t,\left| R\right| \right) \neq \left( 2,16\right) $.
\end{proof}

\textbf{Proof of Proposition \ref{order for GQ}$(iii)$.} If $t=4$,
then there are $16$ points of $S$ not collinear with a given point
$x$. By Proposition \ref{converse}, $|R^{*}\setminus R_{\psi}| \geq
16$. Thus, $|R|>2^{5}$ and so $|R| =2^{6}$. This completes the
proof.

\begin{corollary}
\label{touseinthelastsection}Let $t=4$ and $Q$ be a (2,2)-subGQ of
$S$. Then $|\langle\psi(Q)\rangle|=2^{5}$.
\end{corollary}

\begin{proof}[\bf Proof]
This follows from Proposition \ref{order for GQ}$(iii)$ and the fact
that $P=\langle Q,x\rangle$ for $x\in P\setminus Q$.
\end{proof}

\begin{proposition}
\label{triad}If $t=2$, then $|R|= 2^{4}$ if and only if
$r_{a}r_{b}r_{c}=1$ for every complete 3-arc $\left\{ a,b,c\right\}$
of $S$.
\end{proposition}

\begin{proof}[\bf Proof]
Let $T=\left\{a,b,c\right\} $ be a complete 3-arc of $S$ and $Q$ be
a (2,1)-subGQ of $S$ containing $a$ and $b$. Then $c\notin Q$ and
$P=\langle Q,c\rangle$.

If $r_{a}r_{b}r_{c}=1$, then $r_{c}\in \langle\psi(Q)\rangle$ and
$|R|=|\langle\psi(Q)\rangle|=2^{4}$. Now, assume that $|R|=2^{4}$.
Let $\{x,y\}=\{a,b\}^{\perp}$ in $Q$. Then $x,y\in T^{\perp }$,
since $T$ is a complete 3-arc. Let $z$ be the point in $Q$ such that
$\{ x,y,z\} $ is a 3-arc in $Q$. Then $c\sim z$ and
$r_{z}=(r_{a}r_{x})(r_{b}r_{y})$. Since $ H=\langle r_y:y\in
x^{\perp}\rangle$ is a maximal subgroup of $R$ (\cite{PT}, 4.2.4,
p.68), $|H|=2^{3}$. So $r_{c}=r_{a}r_{b}$ or $r_{a}r_{b}r_{x}$,
since $\psi$ is faithful. If the later holds then $r_{c\ast
z}=r_{y}$, which is not possible because $\psi $ is faithful and
$y\neq c\ast z$. Hence $r_{c}=r_{a}r_{b}$.
\end{proof}

\begin{corollary}
\label{use-1}Assume that $(t,|R|)=(2,2^{4})$. Let
$T=\{a,b,c\}\subset P$ be such that $r_{a}r_{b}r_{c}=1$. Then $T$ is
a line or a complete 3-arc.
\end{corollary}

\begin{proof}[\bf Proof]
Assume that $T$ is not a line. Then, since $\psi$ is faithful, $T$
is a 3-arc. We show that $T$ is complete. Suppose that $T$ is not
complete. Let $\{a,b,d\}$ be the complete 3-arc of $S$ containing
$a$ and $b$. Then $r_{a}r_{b}r_{d}=1$ (Proposition \ref{triad}) and
$c\neq d$. So $r_{c}=r_{d}$, contradicting that $\psi$ is faithful.
\end{proof}

\begin{lemma}
\label{use1}If $S$ contains a 3-arc $T=\{a,b,c\}$ such that
$r_{a}r_{b}r_{c}\in R_{\psi}$, then $(t,|R|)=(2,2^{4})$. In
particular, $T$ is incomplete.
\end{lemma}

\begin{proof}[\bf Proof]
Let $x\in P$ be such that $r_{x}=r_{a}r_{b}r_{c}$. Since $\psi$ is
faithful, $x\notin T$. Let $t=2$. If $T$ is complete, then
$|R|=2^{5}$ (Proposition \ref{triad}) and $x$ is collinear with at
least one point of $T$, say $x\sim a$. Then
$r_{b}r_{c}=r_{x}r_{a}\in R_{\psi}$, a contradiction to Proposition
\ref{converse}. Thus, $T$ is incomplete if $t=2$.

Let $Q_{1}$ be the (2,1)-subGQ of $S$ containing $T$. If $x\in
Q_{1}$, then $\langle \psi(Q_{1})\rangle=\langle
r_{a},r_{b},r_{c},r_{x}\rangle$ would be of order $2^{4}$,
contradicting Proposition \ref{order for GQ}$(i)$. So $x\notin
Q_{1}$ and $t\neq 1$. Let $Q_{2}$ be the (2,2)-subGQ of $S$
generated by $Q_{1}$ and $x$. Then $|\langle
\psi(Q_{2})\rangle|=2^{4}$, and so $t\neq 4$. Thus $t=2$ and
$|R|=|\langle \psi(Q_{2})\rangle|=2^{4}$.
\end{proof}

\begin{lemma}
\label{to use}Let $a,b\in P$ with $a\nsim b$. Set
$A=\{r_{a}r_{x}:x\nsim a\}$ and $B=\{r_{b}r_{x}:x\nsim b\}$. Then
$|A\cap B|=t+2$.
\end{lemma}

\begin{proof}[\bf Proof]
It is enough to prove that $r_{a}r_{x}=r_{b}r_{y}$ for
$r_{a}r_{x}\in A,r_{b}r_{y}\in B$ if and only if either $x=b$ and
$y=a$ holds or there exists a point $c$ such that $\{c,a,y\}$ and
$\{c,b,x\}$ are lines. We need to prove the `only if' part. Since
$\psi$ is faithful, $x\neq b$ if and only if $y\neq a$. Assume that
$x\neq b$ and $y\neq a$. For this, we show that $y\sim a$ and $x\sim
b$. Then $r_{a*y}=r_{a}r_{y}=r_{b}r_{x}=r_{b*x}$. Since $\psi$ is
faithful, it would then follow that $a*y =b*x$ and this would be our
choice of $c$.

First, assume that $(t,|R|)\neq (2,2^{4})$. Since $a\nsim b$,
$r_{a}r_{b}\notin R_{\psi}$ by Proposition \ref{converse}. Since
$r_{x}r_{y}=r_{a}r_{b}$, Proposition \ref{converse} again implies
that $x\nsim y$. Now, $r_{a}r_{b}r_{y}=r_{x}\in R_{\psi}$. By Lemma
\ref{use1}, $\{a,b,y\}$ is not a 3-arc. This implies that $y\sim a$.
By a similar argument, $x\sim b$.

Now, assume that $(t,|R|)= (2,2^{4})$. Suppose that $x\nsim b$.
Then $T=\{a,b,x\}$ is a 3-arc of $S$. By Proposition \ref{use1},
$T$ is incomplete. Let $Q$ be the $(2,1)$-subGQ in $S$ containing
$T$ and let $\{c,d\}=\{a,b\}^{\perp}$ in $Q$. Then
$r_{x}=r_{a}r_{b}r_{c}r_{d}=r_{x}r_{y}r_{c}r_{d}$. So
$r_{y}r_{c}r_{d}=1$. By Corollary \ref{use-1}, $\{c,d,y\}$ is a
complete 3-arc. Since $b\in \{c,d\}^{\perp}$, it follows that
$b\in \{c,d,y\}^{\perp}$, a contradiction to that $b\nsim y$. So
$x\sim b$. A similar argument shows that $y\sim a$.
\end{proof}

\begin{proposition}
\label{expression}Let $K=R^{*}\setminus R_{\psi}$. Each element of
$K$ is of the form $r_{y}r_{z}$ for some $y\nsim z$ in $P$, except
when $(t,|R|)=(2,2^{5})$. In this case, exactly one element, say
$\alpha$, of $K$ can not be expressed in this way. Moreover,
$\alpha =r_{u}r_{v}r_{w}$ for every complete 3-arc $\{u,v,w\} $ of
$S$.
\end{proposition}

\begin{proof}[\bf Proof]
Since $K$ is empty when $(t,|R|)=(2,2^{4})$, we assume that
$(t,|R|)=\left( 1,2^{4}\right) $, $\left( 2,2^{5}\right) $ or
$\left( 4,2^{6}\right) $. Fix $a,b\in P$ with $a\nsim b$. Then
$r_{a}r_{b}\in K$ (Proposition \ref{converse}). Let $A$ and $B$ be
as in Lemma \ref{to use}, and set
\begin{eqnarray*}
C &=&\left\{ r_{a}r_{b}r_{x}:\left\{ a,b,x\right\} \text{is a 3-arc
which is incomplete if }t=2\right\} .
\end{eqnarray*}
By proposition $\ref{converse}$, $A\subseteq K$ and $B\subseteq K$
and by Lemma \ref{use1}, $C\subseteq K$. Each element of $C$
corresponds to a 3-arc which is contained in a (2,1)-subGQ of $S$.
Let $r_{a}r_{b}r_{x}\in C$ and $Q$ be the (2,1)-subGQ of $S$
containing the 3-arc $\left\{ a,b,x\right\} $. If
$\{a,b\}^{\perp}=\{p,q\}$ in $Q$, then $r_{a\ast p}r_{b\ast
q}=r_{x}$ implies that $ r_{a}r_{b}r_{x}=r_{p}r_{q}.$ Thus, every
element of $C$ can be expressed in the required form.

By Proposition \ref{converse}, $A\cap C$ and $B\cap C$ are empty. By
Lemma \ref{to use}, $|A\cap B|= t+2$. Then an easy count shows that
\begin{equation*}
|A\cup B\cup C|=\left\{\begin{array}{ll}
  10t-4&\text{ if }t=1\text{ or }4 \\
  10t-5& \text{ if }t=2
\end{array}.\right.
\end{equation*}
So $K =A\cup B\cup C$ if $t=1$ or $4$, and $K\setminus (A\cup B\cup
C)$ is a singleton if $t=2$. This proves the proposition for $t=1,4$
and tells that if $(t,|R|)=(2,2^{5})$, then at most one element of
$K$ can not be written in the desired form.

Now, let $(t,|R|)=(2,2^{5})$ and $T=\{u,v,w\}$ be a complete 3-arc
of $S$. By Lemma \ref{use1}, $\alpha =r_{u}r_{v}r_{w}\in K$. Suppose
that $ \alpha =r_{x}r_{y}$ for some $x,y\in P$. Then $x\nsim y$ by
Lemma \ref{use1} and $\{x,y\} \cap T=\Phi$ by Proposition
\ref{converse}. Suppose that $x\in T^{\perp}$ and $Q$ be the
$(2,1)$-subGQ of $S$ generated by $\{x,u,v,y\}$. Since $w\notin Q$
and $r_{w}=r_{u}r_{v}r_{x}r_{y}$, it follows that $|R|=2^{4}$, a
contradiction. So, $x\notin T^{\perp}$. Similarly, $y\notin
T^{\perp}$. Thus, each of $x$ and $y$ is collinear with exactly one
point of $T$. Let $x\sim u$. Then $y\nsim x*u$, since $x*u\in
T^{\perp}$ and $\alpha =r_{x}r_{y}$. Let $U$ be the (2,1)-subGQ of
$S$ generated by $\{u,x,y,v\}$. Note that $y\sim u$ in $U$. Let $z$
be the unique point in $U$ such that $\{u,v,z\} $ is a 3-arc of $U$.
Then $r_{z}=r_{x}r_{y}r_{u}r_{v}=r_{w}$. Since $w\neq z$ (in fact,
$w\notin U$), this is a contradiction to the faithfulness of $\psi$.
Thus, $\alpha$ can not be expressed as $r_{x}r_{y}$ for any $x,y$ in
$P$. This, together with the last sentence of the previous
paragraph, implies that $\alpha$ is independent of the complete
3-arc $T$ of $S$.
\end{proof}

\section{Initial Results}

Let $S=(P,L)$ be a slim dense near hexagon and $(R,\psi)$ be a
non-abelian representation of $S$. For $x\in P$ and $y\in
\Gamma_{\leq 2}(x)$, $[r_{x},r_{y}]=1:$ if $d(x,y)=2$, we apply
Proposition \ref{begin} to the restriction of $\psi$ to the quad
$Q(x,y)$. From (\cite{SS}, Theorem 2.9, see Example 2.2 of
\cite{SS}) applied to $S$ , we have

\begin{proposition}
\label{general}\text{ }
\begin{enumerate}
\item[$(i)$] For $x,y\in P$, $[r_{x},r_{y}] \neq 1$ if and only if
$d(x,y)=3$. In this case, $\langle r_{x},r_{y}\rangle$ is a
dihedral group $2_{+}^{1+2}$ of order 8. \item[$(ii)$] $R$ is a
finite 2-group of exponent 4, $|R'|=2$ and $R' =\Phi (R)\subseteq
Z(R)$. \item[$(iii)$] $r_{x}\notin Z(R)$ for each $x\in P$ and
$\psi$ is faithful.
\end{enumerate}
\end{proposition}

We write $R'=\langle \theta\rangle$ throughout. Since $R'$ is of
order two, Lemma \ref{order-calculation} implies

\begin{corollary}
\label{upperbound}$|R|\leq 2^{1+\text{ dim } V(S)}$.
\end{corollary}

\begin{proposition}
\label{centralproduct}$R$ is a central product $E \circ Z(R)$ of
an extraspecial 2-subgroup $E$ of $R$ and $Z(R)$.
\end{proposition}

\begin{proof}[\bf Proof]
We consider $V=R/R'$ as a vector space over $F_{2}$. The map
$f:V\times V\longrightarrow F_{2}$ taking $(xZ,yZ)$ to 0 or 1
accordingly $[x,y] = 1$ or not, is a symplectic bilinear form on
$V$. This is non-degenerate if and only if $R'=Z(R)$. Let $W$ be a
complement in $V$ of the radical of $f$ and $E$ be its inverse image
in $R$. Then $E$ is extraspecial and the proposition follows.
\end{proof}

\begin{corollary}
\label{lowerbound}Let $M$ be an abelian subgroup of $R$ of order
$2^{m}$ intersecting $Z(R)$ trivially. Then $|R|\geq 2^{2m+1}$.
Further, equality holds if and only if $R$ is extraspecial and $M$
is a maximal abelian subgroup of $R$ intersecting $Z(R)$ trivially.
\end{corollary}

The following lemma is useful for us.

\begin{lemma}
\label{useful}Let $x\in P$ and $Y\subseteq \Gamma _{3}(x)$. Then
$[r_{x},\underset{y\in Y}{\Pi }r_{y}] =1$ if and only if $|Y|$ is
even.
\end{lemma}

\begin{proof}[\bf Proof]
Since $R'\subseteq Z(R)$, $[r_{x},\underset{y\in Y}{\Pi }r_{y}]$ is
well-defined (though $\underset{y\in Y}{\Pi }r_{y}$ depends on the
order of multiplication). Let $y,z\in \Gamma _{3}(x)$ be distinct.
The subgraph of $\Gamma(P)$ induced on $\Gamma _{3}(x)$ is connected
(see \cite{BW}, Corollary to Theorem 3, p. 156). Let
$y=y_{0},y_{1},\cdot\cdot\cdot,y_{k}=z$ be a path in $\Gamma
_{3}(x)$. Then $r_{y}r_{z}= {\Pi }r_{y_{i}*y_{i+1}}$ $(0\leq i\leq
k-1)$. Since $d(x,y_{i}*y_{i+1})=2$, $[r_{x},r_{y}r_{z}] =1$. Now,
the result follows from Theorem \ref{general}$(i)$.
\end{proof}

\begin{notation}
For a quad $Q$ in $S$, we denote by $M_{Q}$ the elementary abelian
subgroup of $R$ generated by $\psi(Q)$.
\end{notation}

\begin{proposition}\label{intersectionwithcenter}
Let $Q$ be a quad in $S$ and $M_{Q}\cap Z(R)\neq \{1\}$. Then $Q$ is
of type (2,2), $|M|=2^{5}$ and $M_{Q}\cap
Z(R)=\{1,r_{a}r_{b}r_{c}\}$ for every complete 3-arc $\{a,b,c\}$ of
$S$.
\end{proposition}

\begin{proof}[\bf Proof]
Suppose that $M_{Q}\cap Z(R)\neq\{1\}$ and $1\neq m\in M_{Q}\cap
Z(R)$. Then $m\neq r_{x}$ for each $x\in P$ (Proposition
\ref{general}$(iii)$). If $Q$ is of type (2,1) or (2,4), then by
Proposition \ref{expression}, $m=r_{y}r_{z}$ for some $y,z\in Q,
y\nsim z$. Choose $w\in P\setminus Q$ with $w\sim y.$ Then
$[r_{w},r_{z}]=[r_{w},r_{y}r_{z}]=[ r_{w},m] =1$. But $d(w,z) =3$, a
contradiction to Proposition \ref{general}$(i)$.

So $Q$ is a (2,2)-GQ. Now, $|M_{Q}|\neq 2^{4}$ otherwise
$M_{Q}^{\ast}=\{r_{x}:x\in Q\}$ and $m=r_{x}\in Z(R)$ for some $x\in
Q$, contradicting Proposition \ref{general}$(iii)$. So
$|M_{Q}|=2^{5}$. Now, either $m=r_{u}r_{v}$ for some $u,v\in Q,
u\nsim v$ or $m=r_{a}r_{b}r_{c}$ for every complete 3-arc
$\{a,b,c\}$ of $Q$ (Proposition \ref{expression}). The above
argument again implies that the first possibility does not occur.
\end{proof}

\begin{corollary}\label{intersetion-big-quads}
Let $Q$ and $Q'$ be two disjoint big quads in $S$ of type
$(2,t_{2})$, $t_{2}\neq 2$. Then $M_{Q}\cap M_{Q'}=\{1\}$.
\end{corollary}

\begin{proof}[\bf Proof]
This follows from the proof of Proposition
\ref{intersectionwithcenter} with $Z(R)$ replaced by $M_{Q'}$ and
choosing $w$ in $Q'$.
\end{proof}

\begin{proposition}\label{ovoidal}
Let $Q$ be a quad in $S$ of type $(2,2)$. Then $Q$ is ovoidal if and
only if $|M_{Q}| =2^{5}$ and $M_{Q}\cap Z(R)=\{1 \}$.
\end{proposition}

\begin{proof}[\bf Proof]
First, assume that $Q$ is ovoidal and let $z\in P\setminus Q$ be
such that the pair $(z,Q)$ is ovoidal. Let $\mathcal{O}_{z}=\{
x_{1},\cdot \cdot \cdot ,x_{5}\} $ be as in Theorem \ref{SY1}$(ii)$.
If $|M_{Q}| =2^{4},$ then for the complete 3-arc $\{x_{1},x_{2},y\}$
of $Q$ containing $x_{1}$ and $x_{2}$, $d(y,z)=3$ and
$r_{x_{1}}r_{x_{2}}r_{y}=1$ (Proposition \ref{triad}). But
$[r_{z},r_{y}] =[r_{z},r_{x_{1}}r_{x_{2}}r_{y}] =1$, a contradiction
to Proposition \ref{general}$(i)$. So $|M_{Q}| =2^{5}$. Suppose that
$M_{Q}\cap Z(R)\neq \{1 \}$ and $1\neq m\in M_{Q}\cap Z(R)$. By
Proposition \ref{intersectionwithcenter}, $m=r_{a}r_{b}r_{c}$ for
each complete 3-arc $\{a,b,c\} $ of $Q$. The above argument again
implies that this is not possible. So $M_{Q}\cap Z(R) =\{ 1\}$.

Now, assume that $|M_{Q}|=2^{5}$ and $M_{Q}\cap Z(R)=\{1\}$.
Suppose that $Q$ is classical and let $\{a,b,c\}$ be a complete
3-arc of $Q$. Then, by Proposition \ref{triad},
$r_{a}r_{b}r_{c}\neq 1$. Since $(x,Q)$ is classical for each $x\in
P\setminus Q$, either each of $a,b,c$ is at a distance two from
$x$ or exactly two of them are at a distance three from $x$. In
either case, $[r_{x},r_{a}r_{b}r_{c}]=1$ (see Lemma \ref{useful}).
So $1\neq r_{a}r_{b}r_{c}\in M_{Q}\cap Z(R)$, a contradiction.
\end{proof}

\section{Proof of Theorem \ref{main}}

Let $S=\left(P,L\right)$ be a slim dense near hexagon and let
$(R,\psi)$ be a non-abelian representation of $S$. By Proposition
\ref{general}$(ii)$, $R$ is a finite 2-group of exponent 4. By
Corollary \ref{upperbound}, $|R|\leq 2^{1+dimV(S)}$. For each of the
near hexagons in Theorem \ref{classification-result} except $(vi)$,
we find an elementary abelian subgroup of $R$ of order $2^{\xi}$,
$2\xi =NPdim(S)$, intersecting $Z(R)$ trivially. Then by Corollary
\ref{lowerbound}, $|R|\geq 2^{1+2\xi}$ and $R=2_{+}^{1+2\xi}$ if
equality holds. For the near hexagon $(vi)$ we prove in Subsection
\ref{exceptional-case} that $R=2_{-}^{1+2\xi}$, thus completing the
proof of Theorem \ref{main}.

\subsection{The near hexagons $(vii)$ to $(xi)$}

Let $S=(P,L)$ be one of the near hexagons $(vii)$ to $(xi)$ and $Q$
be a big quad in $S$. Set $M=M_{Q}$. Then, by Proposition
\ref{intersectionwithcenter}, $M\cap Z(R)=\{1\}$ and $|M|=2^{4}$ or
$2^{6}$ according as $Q$ is of type (2,1) or (2,4). If $Q$ is of
type (2,2), then $|M|=2^{4}$ or $2^{5}$. Also, if $|M|=2^{5}$, then
$|M\cap Z(R)|=2$ because $Q$ is classical (Propositions
\ref{intersectionwithcenter} and \ref{ovoidal}). Thus, $R$ has an
elementary abelian subgroup of order $2^{2\xi /2}$ intersecting
$Z(R)$ trivially.

\subsection{The near hexagons $(i)$ to $(v)$}

Let $S=(P,L)$ be one of the near hexagons $(i)$ to $(v)$. Fix $a\in
P$ and $b\in \Gamma _{3}(a)$. Let $l_{1},\cdot \cdot \cdot ,l_{t+1}$
be the lines containing $a$, $x_{i}$ be the point in $l_{i}$ with
$d(b,x_{i}) =2$ and $A= \{ x_{i}:1\leq i\leq t+1\} $. For a subset
$X$ of $A$, we set $T_{X}=\{r_{x}:x\in X \}$, $M_{X}=\langle
T_{X}\rangle$ and $M=\langle r_{b} \rangle M_{X}$. Then $M_{X}$ and
$M$ are elementary abelian 2-subgroups of $R$.

\begin{proposition}
\label{lowerbound1}Let $X$ be a subset of $A$ such that
\begin{enumerate}
\item[$(i)$] $M_{X}\cap Z(R)$=\{1\}, \item[(ii)] $T_{X}$ is
linearly independent.
\end{enumerate}
Then, $|M|=2^{|X| +1}$ and $M\cap Z(R)=\{1\}$. In particular,
$|R|\geq 2^{2|X|+3}$.
\end{proposition}

\begin{proof}[\bf Proof]
By $(ii)$, $2^{|X|}\leq |M| \leq 2^{|X| +1}$. If $|M| =2^{|X|}$,
then $r_{b}$ can be expressed as a product of some of the elements
$r_{x}$, $x\in X$. Since $[r_{a},r_{x}]=1$ for $x\in X$, it follows
that $[r_{a},r_{b}] =1$, a contradiction to Proposition
\ref{general}$(i)$. So $|M| =2^{|X| +1}$. Suppose that $M\cap Z(R)
\neq \{1\}$ and $1\neq z\in M\cap Z(R)$. Let $z=\underset{y\in X\cup
\{ b\} }{\Pi }r_{y}^{i_{y}}$, $i_{y}\in \{0, 1\}$. Since $z\in
Z(R)$, $i_{b}=0$ by the previous argument. Then it follows that
$z\in M_{X}$, a contradiction to $(i)$. So $M\cap Z(R) =\{1\}$.

By Corollary \ref{lowerbound}, $|R|\geq
2^{2(|X|+1)+1}=2^{2|X|+3}$.
\end{proof}

A subset $X$ of $A$ is \textit{good} if $(i)$ and $(ii)$ of
Proposition \ref{lowerbound1} hold. In the rest of this Section, we
find good subsets of $A$ of size $(2\xi -2)/2$, thus completing the
proof of Theorem \ref{main} for the near hexagons $(i)$ to $(v)$.
The next Lemma gives a necessary condition for a subset of $A$ to be
good.

\begin{lemma}\label{[i,j]}
Let $X$ be a subset of $A$ which is not good, $\alpha\in M_{X}\cap
Z(R)$ (possibly $\alpha =1$) and
\begin{equation}\label{expression-eqn}
\alpha =\underset{x_{k}\in X}{\Pi }r_{x_{k}}^{i_{k}}
\end{equation}
where $i_{k}\in \{0,1\}$. Set $B=\left\{k:x_{k}\in X\right\}$,
$B'=\left\{k\in B:i_{k}=1\right\}$ and; for $1\leq i\neq j\leq t+1$,
let $A_{i,j}=\{k\in B':x_{k}\in Q(x_{i},x_{j})\}$. Then

\begin{enumerate}
\item[$(i)$] $|B'|\geq 3$,

\item[$(ii)$] $|B'|$ is even if and only if $|A_{i,j}|$ is even.
\end{enumerate}
\end{lemma}

\begin{proof}[\bf Proof]
$(i)$ $|B'|\geq 2$ because $r_{x_{k}}\notin Z(R)$ for each $k$
(Proposition \ref{general}$(iii)$). If $|B'|= 2$, then
$r_{x}r_{y}=\alpha$ for some pair of distinct $x,y\in X$. Since
$\psi$ is faithful and $r_{x},r_{y}$ are involutions, $\alpha \neq
1$. For the quad $Q=Q(x,y)$, $1\neq \alpha\in M_{Q}\cap Z(R)$. By
Proposition \ref{intersectionwithcenter}, $Q$ is a $(2,2)$-GQ and
$r_{a}r_{b}r_{c}= \alpha$ for each complete 3-arc $\{a,b,c\}$ of
$Q$. In particular, if $\{x,y,w\}$ is the complete 3-arc of $Q$
containing $x$ and $y$, then $r_{x}r_{y}r_{w}= \alpha$. Then it
follows that $r_{w}=1$, a contradiction. So $|B'|\geq 3$.

$(ii)$ Let $w\in Q\left( x_{i},x_{j}\right)$ and $w\nsim a$. For
each $m\in B'_{i,j}=B'\setminus A_{i,j}$, $d(w,x_{m})=3$ because
$x_{m}\sim a$. Now, $[r_{w},\underset{m\in B'_{i,j}}{\Pi
}r_{x_{m}}]=[r_{w},\underset{m\in B'}{\Pi
}r_{x_{m}}]=[r_{w},\alpha]=1$. So $|B'_{i,j}|$ is even by Lemma
\ref{useful}. This implies that $(ii)$ holds.
\end{proof}

In what follows, for any subset $X$ of $A$ which is not good, $B'$
is defined relative to an expression as in (\ref{expression-eqn})
for an arbitrary but fixed element of $M_{X}\cap Z(R)$. Any quad $Q$
in $S$ containing the point $a$ is determined by any two distinct
points $x_{i}$ and $x_{j}$ of $A$ that are contained in $Q$. In that
case we sometime denote by $A_{Q}$ the set $A_{i,j}$ defined in
Lemma \ref{[i,j]}.

\subsubsection{The near hexagon $(i)$}

There are 7 quads in $S$ containing the point $x_{1}\in A$. This
partitions the 14 points ($\neq x_{1}$) of $A$ , say
\begin{equation*}
\{x_{2},x_{3}\} \cup \{x_{4},x_{5}\} \cup \{x_{6},x_{7}\} \cup
\{x_{8},x_{9}\} \cup \{x_{10},x_{11}\} \cup \{x_{12},x_{13}\} \cup
\{x_{14},x_{15}\} .
\end{equation*}
Consider the quad $Q(x_{10},x_{12})$. We may assume that
$Q(x_{10},x_{12})\cap A=\{x_{10},x_{12},x_{15}\}$. We show that
$X=\{x_{2},x_{3},x_{4},x_{5},x_{6},x_{7},x_{8},x_{10},x_{12},x_{14}\}$
is a good subset of $A$.

Assume otherwise. Let $C_{1}=\{8,10,12,14\}$ and $C_{2}=B\setminus
C_{1}$. For $k\in C_{1}$, $Q(x_{1},x_{k})\cap
A=\{x_{1},x_{k},x_{k+1}\}$. So $A_{1,k}\subseteq\{k\}$. By Lemma
\ref{[i,j]}$(ii)$, either $C_{1}\subseteq B'$ or $C_{1}\cap B'$ is
empty. Now, $C_{1}\nsubseteq B'$ because, otherwise, $A_{1,14}=\{14
\}$ and $A_{10,12}= \{10,12\}$ and, by Lemma \ref{[i,j]}$(ii)$,
$|B'|$ would be both odd and even.

Suppose that $C_{1}\cap B'$ is empty. Then $B'\subseteq C_{2}$.
Since $A_{1,8}$ is empty, $|B'|$ is even. Choose $j\in B'$ (see
Lemma \ref{[i,j]}$(i)$). Observe that there exists $k\in
\{8,\cdot\cdot\cdot,15\}$ such that $Q(x_{j},x_{k})\cap
\{x_{i}:i\in C_{2}\}=\{x_{j}\}$. Then $A_{j,k}=\{j\}$ and $|B'|$
is odd also, a contradiction. So, $X$ is good and $|X|=10$.

\subsubsection{The near hexagon $(ii)$}

Let $X=\{x_{i}:1\leq i\leq 11\}$. Then $X$ is a good subset of $A$.
Otherwise, for some $i,j\in B'$ with $i\neq j$ (see Lemma
\ref{[i,j]}$(i)$), $A_{i,j}=\{i,j\}$ and $A_{i,12}=\{i\}$ and, by
Lemma \ref{[i,j]}$(ii)$, $|B'|$ would be both even and odd.

\subsubsection{The near hexagon $(iii)$}

Let $Q_{1},\cdot \cdot \cdot ,Q_{5}$ be the five (big) quads in $S$
containing $x_{1}$ and $a$. Let
\begin{center}
\begin{tabular}{l}
$Q_{1}\cap A=\{x_{1},x_{2},x_{3},x_{4},x_{5}\},$ \\
$Q_{2}\cap A=\{x_{1},x_{6},x_{7},x_{8},x_{9}\},$ \\
$Q_{3}\cap A=\{x_{1},x_{10},x_{11},x_{12},x_{13}\},$ \\
$Q_{4}\cap A=\{x_{1},x_{14},x_{15},x_{16},x_{17}\},$ \\
$Q_{5}\cap A=\{x_{1},x_{18},x_{19},x_{20},x_{21}\}.$
\end{tabular}
\end{center}
We show that
$X=\{x_{2},x_{3},x_{4},x_{5},x_{6},x_{7},x_{8},x_{10},x_{14}\}$ is a
good subset of $A$. Assume otherwise. Since $Q_{5}\cap X$ is empty,
$A_{Q_{5}}$ is empty and, by Lemma \ref{[i,j]}$(ii)$, $|B'|$ and
$|A_{Q}|$ are even for each quad $Q$ in $S$ containing $a$. Since
$A_{Q_{3}}\subseteq \{10\}$ and $|A_{Q_{3}}|$ is even, $10\notin
A_{Q_{3}}$ and so, $10\notin B'$. This argument with $Q_{3}$
replaced by $Q_{4}$ shows that $14\notin B'$. Since
$A_{Q_{2}}\subseteq \{6,7,8\}$ and $|A_{Q_{2}}|$ is even, $j\notin
B'$ for some $j\in \{6,7,8\}$. Since $|B'|\geq 3$ (Lemma
\ref{[i,j]}$(i)$), $k\in B'$ for some $k\in \{2,3,4,5\}$. Then,
$A_{j,k}=\{k\}$, contradicting that $|A_{j,k}|$ is even. So $X$ is
good and $|X|=9$.

\subsubsection{The near hexagon $(iv)$}

Let $Q_{1},\cdot \cdot \cdot ,Q_{6}$ be the six big quads in $S$
containing the point $a$. Any two of these big quads meet in a line
through $a$ and any three of them meet only at $\{a\}$. Let
\begin{center}
\begin{tabular}{l}
$Q_{1}\cap A=\{x_{1},x_{2},x_{3},x_{4},x_{5}\},$\\
$Q_{2}\cap A=\{x_{1},x_{6},x_{7},x_{8},x_{9}\},$\\
$Q_{3}\cap A=\{x_{2},x_{6},x_{10},x_{11},x_{12}\},$\\
$Q_{4}\cap A=\{x_{3},x_{7},x_{10},x_{13},x_{14}\},$\\
$Q_{5}\cap A=\{x_{4},x_{8},x_{11},x_{13},x_{15}\},$\\
$Q_{6}\cap A=\{x_{5},x_{9},x_{12},x_{14},x_{15}\}.$
\end{tabular}
\end{center}
We show that
$X=\{x_{1},x_{2},x_{3},x_{4},x_{6},x_{7},x_{8},x_{10},x_{11}\}$ is a
good subset of $A$. Assume otherwise. Since $Q_{6}\cap X$ is empty,
$A_{Q_{6}}$ is empty and, by Lemma \ref{[i,j]}$(ii)$, $|B'|$ and
$|A_{Q}|$ are even for every quad $Q$ in $S$ containing $a$. We
first verify that for
$$(i,j,k)\in
\{(1,11,14),(1,12,13),(2,9,13),(3,6,15),(4,6,14),(5,6,13)\},$$
$Q(x_{i},x_{j})$ is of type (2,2) and $Q(x_{i},x_{j})\cap
A=\{x_{i},x_{j},x_{k}\}$. Since $A_{1,12}\subseteq \{1\}$ and
$|A_{1,12}|$ is even, it follows that $1\notin B'$. Similarly,
considering $A_{2,9}$ and $A_{5,6}$, we conclude that $2\notin B'$
and $6\notin B'$. Since $6\notin B'$, considering $A_{3,6}$ and
$A_{4,6}$, we conclude that $3\notin B'$ and $4\notin B'$. Since
$|B'|\geq 3$ is even, it follows that $B'=\{7,8,10,11\}$ and so
$A_{1,11}=\{11\}$, contradicting that $|A_{1,11}|$ is even. So $X$
is good and $|X|=9$.

\subsubsection{The near hexagon $(v)$}

Let $Q_{1},Q_{2},Q_{3}$ be the three big quads containing $a$. There
intersection is $\{a\}$ and any two of these big quads meet in a
line through $a$. We may assume that
\begin{center}
\begin{tabular}{l}
$Q_{1}\cap A=\{x_{1},x_{2},x_{3},x_{4},x_{5}\},$\\
$Q_{2}\cap A=\{x_{1},x_{6},x_{7},x_{8},x_{9}\},$\\
$Q_{3}\cap A=\{x_{2},x_{6},x_{10},x_{11},x_{12}\}.$
\end{tabular}
\end{center}
We show that
$X=\{x_{1},x_{2},x_{3},x_{4},x_{6},x_{7},x_{8},x_{10},x_{11}\}$ is
good subset of $A$. Assume otherwise. We note that the quads
$Q(x_{r},x_{k})$ are of type (2,2) in the following cases:
$$r=1\text{ and }k\in\{10,11,12\}; r=2\text{ and }k\in\{7,8,9\};
r=6\text{ and }k\in\{3,4,5\}.$$ Now, $A_{r,s}\subseteq \{r\}$ for
$(r,s)\in\{(1,12),(2,9),(6,5)\}$ because $x_{s}\notin X$.
Considering $A_{1,12}$, we conclude that $10,11\notin B'$ in view of
the following: $A_{1,12}\subseteq \{1\}$, $A_{1,k}\subseteq\{1,k\}$
for $k\in\{10,11\}$ and the parity of $|B'|$ and $|A_{1,j}|$ are the
same for all $j\neq 1$. Similarly, considering $A_{2,9}$
(respectively, $A_{6,5}$) we conclude that $7,8\notin B'$
(respectively, $3,4\notin B'$). Since $|B'|\geq 3$, it follows that
$B'=\{1,2,6\}$. But $A_{5,9}$ is empty because
$\{x_{5},x_{9},x_{12}\}\cap X$ and $\{10,11\}\cap B'$ are empty. So
$|B'|$ is even (Lemma \ref{[i,j]}$(ii)$), a contradiction. So $X$ is
good and $|X|=9$.

\subsection{The near hexagon $(vi)$}\label{exceptional-case}

We consider this case separately because the technique of the
previous section only yields $|R|\geq 2^{17}$ in this case.

Let $S=\left( P,L\right) $ be a slim dense near hexagon and $Y$ be
a proper subspace of $S$ isomorphic to the near hexagon $(vii)$.
Big quads in $Y$ (as well as in $S$) are of type (2,4). There are
three pair-wise disjoint big quads in $Y$ and any two of them
generate $Y$. Fix two disjoint big quads $Q_{1}$ and $Q_{2}$ in
$Y$. Let $(R,\psi)$ be a non-abelian representation of $S$. Set
$M=\langle \psi (Y)\rangle$ and $M_{i}=M_{Q_{i}}$ for $i=1,2$.
Then $|M_{i}| =2^{6}$ (Proposition \ref{order for GQ}$(iii)$),
$M_{i}\cap Z(R) =\{1 \}$ (Proposition
\ref{intersectionwithcenter}), $ M_{1}\cap M_{2}=\{1 \}$
(Corollary \ref{intersetion-big-quads}) and $M=2_{+}^{1+12}$ with
$M=M_{1}M_{2}R^{\prime}$ (Theorem \ref{main} for the the near
hexagon $(vii)$). Clearly, $R=M\circ N$, where $N=C_{R}(M)$.

Let $\{i,j\} =\{1,2\}$. For $x\in P\setminus Y$, we denote by
$x^{j}$ the unique point in $Q_{j}$ at distance 1 from $x$. For
$y\in Q_{i}$, let $z_{y}$ denote the unique point in $Q_{j}$ at
distance 1 from $y$.

\begin{proposition}
\label{done2(vi)}For each $x\in P\setminus Y$, $r_{x}$ has a unique
decomposition as $r_{x}=m_{1}^{x}m_{2}^{x}n_{x}$, where
$m_{j}^{x}=r_{z_{x^{i}}}\in M_{j}$ and $n_{x}\in N$ is an involution
not in $Z(R)$. In particular, $r_{x}\notin M$.
\end{proposition}

\begin{proof}[\bf Proof]
We can write $r_{x}=m_{1}^{x}m_{2}^{x}n_{x}$ for some $m_{1}^{x}\in M_{1},$ $%
m_{2}^{x}\in M_{2}$ and $n_{x}\in N$. Set $H_{j}=\langle r_{w}: w\in
Q_{j}\cap x^{j\perp }\rangle\leq M_{j}$. Then $H_{j}$ is a maximal
subgroup of $M_{j}$ (\cite{PT}, 4.2.4, p.68) and $r_{x}\in
C_{R}(H_{1})\cap C_{R}(H_{2})$. For all $h\in H_{j},$
\begin{equation*}
\left[ m_{i}^{x},h\right]=\left[
m_{1}^{x}m_{2}^{x}n_{x},h\right]=\left[ r_{x},h\right] =1.
\end{equation*}
So $ m_{i}^{x}\in C_{M_{i}}\left( H_{j}\right)$. Note that
$C_{M_{i}}\left(H_{j}\right) =\langle r_{z_{x^{j}}}\rangle$, a
subgroup of order 2. If $m_{i}^{x}=1$, then $r_{x}=m_{j}^{x}n_{x}$
commutes with every element of $ M_{j}$. In particular, $\left[
r_{x},r_{y}\right] =1$ for every $y\in Q_{j}\cap \Gamma _{3}\left(
x\right)$, a contradiction to Theorem \ref{general}$(i)$. So
$m_{i}^{x}=r_{z_{x^{j}}}$. Now $\left[ m_{1}^{x},m_{2}^{x} \right]
=1$, since $d\left( z_{x^{ 1}},z_{x^{ 2}}\right) =2$ (Proposition
\ref{done1(vi)}). Since $r_{x}^{2}=1$, $n_{x}^{2}=1$.

We show that $n_{x}\neq 1$ and $n_{x}\notin Z(R)$. The quad $Q =Q(
x^{1},x^{2})$ is of type (2,2) or (2,4) because $x^{1}$ and $x^{2}$
have at least three common neighbours $x,z_{x^{1}}$ and $z_{x^{2}}$.
Let $U$ be the $\left( 2,2\right) $-GQ in $Q$ generated by $ \left\{
x^{1},x^{2},x,z_{x^{1}},z_{x^{2}}\right\}$. If $Q$ is of type (2,4),
then $\langle \psi (U) \rangle$ is of order $2^{5}$ (Corollary
\ref{touseinthelastsection}). If $Q$ is of type (2,2), then $U=Q$ is
ovoidal because it is not a big quad. So $\langle \psi (U) \rangle$
is of order $2^{5}$ (Propositions \ref{ovoidal}). Therefore,
$r_{a}r_{b}r_{c}\neq 1$ for every complete 3-arc $\left\{
a,b,c\right\} $ of $U$ (Proposition \ref{triad}). In particular,
$n_{x}=r_{x}r_{z_{x^{1}}}r_{z_{x^{2}}}\neq 1$ for the complete 3-arc
$\left\{ x,z_{x^{1}},z_{x^{2}}\right\}$ of $U$. Now, applying
Proposition \ref{intersectionwithcenter} (respectively, Proposition
\ref{ovoidal}) when $Q$ is of type (2,4) (respectively, of type
(2,2)), we conclude that $n_{x} \notin Z(R)$.
\end{proof}

\begin{proposition}
\label{done4(vi)}Let $Q$ be a big quad in $S$ disjoint from $Y$
and $x,y\in Q$. Then:
\begin{enumerate}
\item[$(i)$] $[n_{x},n_{y}]=1$ if and only if $x=y$ or $x\sim y$;

\item[$(ii)$] There is a unique line $l_{x}=\left\{ x,y,x\ast
y\right\} $ in $Q$ containing $x$ such that $ n_{x\ast
y}=n_{x}n_{y}$. For any other line $l=\left\{ x,z,x\ast z\right\} $
in $Q$, $n_{x\ast z}=n_{x}n_{z}\theta$.
\end{enumerate}
\end{proposition}

\begin{proof}[\bf Proof]
$(i)$ Let $x\sim y$. By Corollary \ref{done(vi)-corollary} and
Proposition \ref{done2(vi)},
$[m_{2}^{x},m_{1}^{y}]=[m_{1}^{x},m_{2}^{y}]=1$ or $\theta$ . Then
$[n_{x},n_{y}]=[m_{1}^{x}m_{2}^{x}n_{x},m_{1}^{y}m_{2}^{y}n_{y}]=[r_{x},r_{y}]
=1$.

Now, assume that $x\nsim y$. By Propositions
\ref{done(vi)-proposition} and \ref{done2(vi)},
$([m_{1}^{x},m_{2}^{y}], [m_{2}^{x},m_{1}^{y}])=(1,\theta)$ or
$(\theta,1)$. Since $[r_{x},r_{y}] =1$, it follows that
$[n_{x},n_{y}]=\theta\neq 1$.

$(ii)$ Let $x\in Q$ and $l_{x}$ be the line in $Q$ containing $x$
which corresponds to the line $x^{j}z_{x^{i}}$ in $Q_{j}$. This is
possible by Lemma \ref{isomorphism}. For $u,v\in l_{x}$, $d(
z_{u^{j}},z_{v^{i}}) \leq 2$ (Corollary \ref{done(vi)-corollary}).
So $[m_{i}^{u},m_{j}^{v}] =1$. Then $r_{u\ast v}=(
m_{1}^{u}m_{1}^{v})( m_{2}^{u}m_{2}^{v}) (n_{u}n_{v})$. So $n_{u\ast
v}=n_{u}n_{v}.$ Let $l$ be a line ($\neq l_{x}$) in $Q$ containing
$x$. For $y\neq w$ in $l$, $[m_{2}^{y},m_{1}^{w}] =\theta$ because
$d(z_{y^{1}},z_{w^{2}})=3$ (Corollary \ref{done(vi)-corollary}). So
$r_{y\ast w}=\left( m_{1}^{y}m_{2}^{y}n_{y}\right) \left(
m_{1}^{w}m_{2}^{w}n_{w}\right) =\left( m_{1}^{y}m_{1}^{w}\right)
\left( m_{2}^{y}m_{2}^{w}\right) n_{y}n_{w}\theta$ and $n_{y\ast
w}=n_{y}n_{w}\theta$.
\end{proof}

\begin{corollary}
\label{done5(vi)}Let $Q$ be as in Proposition \ref{done4(vi)} and
$I_{2}(N)$ be the set of involutions in $N$. Define $\delta$ from
$Q$ to $I_{2}(N)$ by $\delta (x)=n_{x}$. Then
\begin{enumerate}
\item[$(i)$] $[\delta (x),\delta (y)] =1$ if and only if $x=y$ or
$x\sim y$.

\item[$(ii)$] $\delta $ is one-one.

\item[$(iii)$] There exists a spread $T$ in $Q$ such that for
$x,y\in Q$ with $x\sim y$,
\begin{equation*}
\delta (x\ast y)=\left\{\begin{array}{ll}
  \delta (x) \delta (y) & \text{ if }xy\in T \\
  \delta (x) \delta (y)\theta & \text{ if }xy\notin T
\end{array}.\right.
\end{equation*}
\end{enumerate}
\end{corollary}

\begin{proof}[\bf Proof]
$(i)$ and $(iii)$ follows from Proposition \ref{done4(vi)}. We now
prove $(ii)$. Let $\delta (x)=\delta(y)$ for $x,y\in Q$. By $(i)$,
$x=y$ or $x\sim y$. If $x\sim y$, then $r_{x\ast
y}=r_{x}r_{y}=(m_{1}^{x}m_{1}^{y})( m_{2}^{x}m_{2}^{y}) \alpha \in
M$, where $\alpha =[m_{2}^{x},m_{1}^{y}] \in R^{\prime }.$ But this
is not possible as $x\ast y\notin Y$ (Proposition \ref{done2(vi)}).
So $x=y$.
\end{proof}

Now, let $S=(P,L)$ be the near hexagon $(vi)$. Then big quads in $S$
are of type (2,4). We refer to (\cite{BCHW}, p.363) for the
description of the corresponding Fischer Space on the set of 18 big
quads in $S$. This set partitions into two families $F_{1}$ and
$F_{2}$ of size $9$ each such that each $F_{i}$ defines a partition
of the point set $P$ of $S$. Let $U_{i}$, $i=1,2$, be the partial
linear space whose points are the big quads of $F_{i}$, two distinct
big quads considered to be collinear if they are disjoint. If
$Q_{1}$ and $Q_{2}$ are collinear in $U_{i}$, then the line
containing them is $\{Q_{1}, Q_{2}, Q_{1}\ast Q_{2}\}$, where
$Q_{1}\ast Q_{2}$ is defined as in Lemma \ref{isomorphism}. Then
$U_{i}$ is an affine plane of order 3.

Consider the family $F_{1}$. Fix a line $\{Q_{1}, Q_{2}, Q_{1}\ast
Q_{2}\}$ in $U_{1}$ and set $Y=Q_{1}\cup Q_{2}\cup Q_{1}\ast Q_{2}$.
Then $Y$ is a subspace of $S$ isomorphic to the near hexagon
$(vii)$. Fix a big quad $Q$ in $U_{1}$ disjoint from $Y$. Let the
subgroups $M$ and $N$ of $R$ be as in the beginning of this
subsection. Then $|N|\leq 2^{7}$ because $|R|\leq
2^{1+dimV(S)}=2^{19}$. We show that $N=2_{-}^{1+6}$. This would
prove Theorem \ref{main} in this case.

Let $\{a_{1},a_{2},b_{1},b_{2}\}$ be a quadrangle in $Q$, where
$a_{1}\nsim a_{2}$ and $b_{1}\nsim b_{2}$. Let $\delta$ be as in
Corollary \ref{done5(vi)}. Then the subgroup $\langle
\delta(a_{1}),\delta(a_{2}),\delta(b_{1}),\delta(b_{2})\rangle $ of
$R$ is isomorphic to $H=\langle \delta(a_{1}),\delta(a_{2})\rangle
\circ \langle \delta(b_{1}),\delta(b_{2})\rangle$. We write $N=H
\circ K$ where $K=C_{N}(H)$. Then $|K|\leq 2^{3}$. There are three
more neighbours, say $w_{1},w_{2},w_{3}$, of $a_{1}$ and $a_{2}$ in
$Q$ different from $b_{1}$ and $b_{2}$. We can write
$$\delta(w_{i})=
\delta(a_{1})^{i_{1}}\delta(a_{2})^{i_{2}}\delta(b_{1})^{j_{1}}\delta(b_{2})^{j_{2}}k_{i}$$
for some $k_{i}\in K$, where $i_{1},i_{2},j_{1},j_{2}\in \{0,1\}$.
By Corollary \ref{done5(vi)}$(i)$,
$[\delta(w_{i}),\delta(a_{r})]=1\neq
[\delta(w_{i}),\delta(b_{r})]$ for $r=1,2$. This implies that
$i_{1}=i_{2}=0$ and $j_{1}=j_{2}=1$; that is, $\delta
(w_{i})=\delta(b_{1})\delta(b_{2})k_{i}$. In particular, $k_{i}$
is of order 4. Since $[\delta(w_{i}),\delta(w_{j}]\neq 1$ for
$i\neq j$, it follows that $[k_{i},k_{j}]\neq 1$. Thus, $K$ is
non-abelian and is of order 8 and $k_{1},k_{2}$ and $k_{3}$ are
three pair-wise distinct elements of order 4 in $K$. So $K$ is
isomorphic to $Q_{8}$ and $N=2_{-}^{1+6}$.

\end{document}